\documentclass[12pt]%
{amsart}
\usepackage{verbatim}
\usepackage{t1enc}
\usepackage[latin2]{inputenc}

\usepackage{amssymb}
\usepackage{enumerate}
\newif\ifdeveloping

\newtheorem{theorem}{Theorem}

\newtheorem{claim}{Claim}[theorem]
\newtheorem{problem}[theorem]{Problem}

\newtheorem{qtheorem}{Theorem}

\theoremstyle{definition}

\theoremstyle{remark}
  
{}

\ifdeveloping
\usepackage[notref,notcite]{showkeys}
\fi

\newcommand{\prtime}{{\count0=\time\divide\count0 by 60
\count1=-\count0\multiply\count1 by 60
\advance\count1 by \time
\the\count0:\the\count1}
}

\def\myheads#1;#2;{
\pagestyle{myheadings}
\markboth{{\sc\hfill #1\hfill\protect\makebox[0cm][r]{\rm\today; \prtime}}}
{{\sc\protect\makebox[0cm][l]{\rm\today;\ \prtime}\hfill #2\hfill}}
\thispagestyle{myheadings}
}

\newcommand{\dcal}{{\mathcal D}}

\newcommand{\ncal}{{\mathcal N}}

\newcommand{\vcal}{{\mathcal V}}

\newcommand{\setm}{\setminus}
\newcommand{\empt}{\emptyset}
\newcommand{\subs}{\subset}
\newcommand{\sups}{\supset}

\newcommand{\dom}{\operatorname{dom}}
\newcommand{\nden}{\operatorname{\ncal}}

\def\<{\left\langle}
\def\>{\right\rangle}

\def\br#1;#2;{\bigl[ {#1} \bigr]^ {#2} }
\newcommand{\fn}{\operatorname{Fn}}

\theoremstyle{plain}

\newcommand{\eres}{extraresolvable}

\newcommand{\nwd}{\operatorname{nwd}}
\newcommand{\nwds}{\operatorname{\ncal}}

\newcommand{\bbb}{\mathbb B}
\newcommand{\cbb}{\mathbb C}
\newcommand{\dbb}{\mathbb{D}}

\newcommand{\mosaic}[1]{$#1$-mosaic}

\newcommand{\cb}{\mathbb C}
\newcommand{\db}{\mathbb D}
\newcommand{\eb}{\mathbb E}
\newcommand{\ap}{a^+}
\newcommand{\am}{a^-}

\begin{document}

\author[I. Juh\'asz]{Istv\'an Juh\'asz}
\address{Alfr{\'e}d R{\'e}nyi Institute of Mathematics}
\email{juhasz@renyi.hu}

\author[S. Shelah]{
Saharon Shelah}
\address{Hebrew University}
\email{shelah@math.huji.ac.il}

\author[L. Soukup]{
Lajos Soukup }
\address{
Alfr{\'e}d R{\'e}nyi Institute of Mathematics }
\email{soukup@renyi.hu}

\title{Resolvability vs. almost resolvability}
\thanks{The first and third author was supported by OTKA grant no. 61600}

\maketitle
\begin{abstract}
A space $X$ is {\em $\kappa$-resolvable} (resp. {\em almost
$\kappa$-resolvable}) if it contains $\kappa$ dense sets that are
pairwise disjoint (resp. almost disjoint over the ideal of nowhere
dense subsets of $X$).

Answering a problem raised by Juh\'asz, Soukup, and Szentmikl\'ossy,
and improving a consistency result of Comfort and Hu, we prove, in
ZFC,  that for every infinite cardinal ${\kappa}$ there is an almost
$2^{\kappa}$-resolvable but not ${\omega}_1$-resolvable  space of
dispersion character ${\kappa}$.
\end{abstract}

A space $X$ is said to be {\em $\kappa$-resolvable} if it contains
$\kappa$ dense sets that are pairwise disjoint. $X$ is called {\em
maximally resolvable} iff it is $\Delta(X)$-resolvable, where
$\Delta(X) = \min\{ |G| : G \ne \emptyset \mbox{ open}\}$ is the
{\em dispersion character} of $X$.

V. Malychin, in \cite{Ma3},  was the first to suggest studying
families of dense sets of a space $X$ that, rather than disjoint,
are merely {\em almost disjoint} with respect to the ideal
$\nden(X)$, where $\nden(X)$ denotes the family of all nowhere dense
subsets of the space $X$. He called a space $X$ {\em\eres} if it has
$\operatorname{\Delta}(X)^+$ many dense sets such that any two of
them have nowhere dense intersection. This idea was generalized in
\cite{d-forced}, where the natural notion of {\em almost
$\kappa$-resolvability} was introduced: A space $X$ is called {\em
almost $\kappa$-resolvable} if it contains $\kappa$ dense sets that
are pairwise almost disjoint over the ideal $\nden(X)$ of nowhere
dense subsets of $X$. (Actually, this concept was given a different
name in \cite{d-forced}, namely: ``$\kappa$-extraresolvable '', but
we think the terminology given here is much better.)

Note that this makes good sense for $\kappa \le
\operatorname{\Delta}(X)$ as well. But while ``almost
${\omega}$-resolvable'' is clearly equivalent to
``${\omega}$-resolvable'', the analogous question for higher
cardinals remained open. In particular, the following natural
problem was formulated in \cite{d-forced}:
\begin{problem}\label{pr:ex}
Let $X$ be an extraresolvable ($T_2$, $T_3,$ or Tychonov) space with
$\operatorname{\Delta}(X) \ge \omega_1$. Is $X$ then $\omega_1$-
resolvable?
\end{problem}
(The assumption $\operatorname{\Delta}(X) \ge \omega_1$ is clearly
necessary to make this problem non-trivial.)

Comfort and Hu, see {\cite[Corollary 3.6]{CoHu2}}, gave a negative
answer to this problem, assuming the failure of the continuum
hypothesis, CH. More precisely they got the following result:
\begin{qtheorem}
If $\kappa$ is an infinite cardinal such that $GCH$ first fails at
${\kappa}$ then there is a 0-dimensional $T_2$ space $X$ with
$|X|=\Delta(X)={\kappa}^+$ such that $X$ is ${\kappa}$-resolvable,
extraresolvable but not ${\kappa}^+$-resolvable, hence not maximally
resolvable and if $\kappa = \omega$ then not $\omega_1$- resolvable.
\end{qtheorem}

Our aim in this note is to give the following ``final'' answer to
the above problem, in ZFC.

\begin{theorem}\label{tm:main1}
For every cardinal ${\kappa}$ there is a 0-dimensional $T_2$ space
of dispersion character ${\kappa}$ that is extraresolvable  but not
${\omega}_1$-resolvable.
\end{theorem}

We shall actually prove a bit more. Note that no space $X$ can be
almost ${(2^{\Delta(X)})}^+$-resolvable, moreover ``almost
$2^{\Delta(X)}$-resolvable'' can be strictly stronger than
``extraresolvable $\equiv \mbox{ almost }
{\Delta(X)}^+$-resolvable''.

\begin{theorem}\label{tm:main2}
For every cardinal ${\kappa}$ there is an almost
$2^{\kappa}$-resolvable (and so extraresolvable) but not
${\omega}_1$-resolvable 0-dimensional $T_2$ space of cardinality and
dispersion character ${\kappa}$. In fact, our example is a
$\kappa$-dense subspace of the Cantor cube of weight $2^\kappa$.
\end{theorem}

To prove this theorem we shall make use of the method of
constructing $\dcal$-forced spaces that was introduced in
\cite{d-forced}. Therefore, we first recall some definitions and
results from \cite{d-forced}.


Let $\dcal$ be a family of dense subsets of a space $X$.  A subset
$M\subs X$ is called a {\em \mosaic {\dcal}} iff
there is a maximal disjoint family
$\vcal$ of open subsets of $X$ and for each $V\in\vcal$
there is $D_V\in \dcal$
such that
\begin{displaymath}
M=\cup\{V\cap D_V:V\in\vcal\}.
\end{displaymath}
Clearly, every \mosaic {\dcal} is dense. We say  that the space $X$
(or its topology)
 is  {\em ${\dcal}$-forced} iff
every dense subset of $X$ includes a $\dcal$-mosaic.

Let $S$ be any set and $\bbb=\left\{\<B_{\zeta}^0,B_{\zeta}^1\> :
{\zeta}<{\mu}\right\}$ be a family of 2-partitions of $S$. We denote
by $\tau_\bbb$ the (obviously zero-dimensional) topology on $S$
generated by the subbase $\{B_{\zeta}^i: \zeta < \mu,\ i < 2\}$,
moreover we set $X_{\bbb}=\<{S},{\tau}_{\bbb}\>$.

Given a cardinal ${\kappa}$, we have $\Delta(X_{\bbb})\ge {\kappa}$
iff $\bbb$ is {\em ${\kappa}$-independent}, i.e.,
\begin{displaymath}
\bbb[\varepsilon]\stackrel{def}=
\bigcap\{B^{{\varepsilon}({\zeta})}_{\zeta}:{\zeta}\in\dom
{\varepsilon} \}
\end{displaymath}
has cardinality at least ${\kappa}$ whenever ${\varepsilon}\in
Fn({\mu},2)$.

Note that $X_{\bbb}$ is Hausdorff iff $\bbb$ is {\em separating},
i.e. for each pair $\{{\alpha},{\beta}\}\in \br S ;2;$ there are
${\zeta}<{\mu}$ and $i<2$ such that ${\alpha}\in B^i_{\zeta}$ and
${\beta}\in B^{1-i}_{\zeta}$.

A set $D\subs X$ is said to be {\em ${\kappa}$-dense} in the space
$X$ iff $|D\cap U|\ge {\kappa}$ for each nonempty open set $U\subs
X$. Thus $D$ is dense iff it is 1-dense. Also, it is obvious that
the existence of a $\kappa$-dense set in $X$ implies
$\operatorname{\Delta}(X) \ge \kappa$.


\begin{qtheorem}[{\cite[Main Theorem 3.3]{d-forced}}]\label{tm:ld}
\label{tm:main} Assume that ${\kappa}$ is an infinite cardinal and
we are given $\bbb= \bigl\{\<B_{\xi}^0,B_{\xi}^1\> :
{\xi}<2^{\kappa} \bigr\}$, a ${\kappa}$-independent family of
$\,2$-partitions of ${\kappa}$,  moreover a non-empty family $\dcal$
of ${\kappa}$-dense subsets of the space $X_{\bbb}$. Then there is a
separating ${\kappa}$-independent family
$\cbb=\{\<C_{\xi}^0,C_{\xi}^1\>:{\xi}< 2^{\kappa}\}$ of
$\,2$-partitions of ${\kappa}$ such that
\begin{enumerate}[(1)]
\item \label{dense}
every $D\in \dcal$ is also ${\kappa}$-dense in $X_{\cbb}$ (and so
$\operatorname{\Delta}(X_\cbb)={\kappa}$),
\item \label{forced} $X_{\cbb}$ is $\dcal$-forced.

\end{enumerate}
\end{qtheorem}
Actually, the space $X_{\cbb}$ has other interesting properties as
well but we shall note make use of those here. We are now ready to
prove our promised result.
\begin{proof}[Proof of Theorem \ref{tm:main2}]
Let ${\kappa}$ be an arbitrary infinite cardinal. It is well-known,
see e. g.  \cite[Fact 3.2]{d-forced}, that we can find two disjoint
families $\mathbb B=\bigl\{\<B_i^0,B_i^1\>:i<2^{\kappa} \bigr\}$ and
$\mathbb D=\bigl\{\<D_i^1,D_i^1\>:i<2^{\kappa}\bigr\}$ of
$2$-partitions of ${\kappa}$ such that their union $\bbb\cup\dbb$ is
${\kappa}$-independent, that is, for any ${\eta},{\varepsilon}\in
\fn (2^{\kappa},2)$ we have
\begin{displaymath}
\bigl|\ \mathbb D[{\eta}]\cap \mathbb B[{\varepsilon}]\
\bigr|={\kappa}.
\end{displaymath}
In other words, this means that
\begin{displaymath}
\dcal=\{\dbb[{\eta}]:{\eta}\in \fn (2^{\kappa},2)\}
\end{displaymath}
is a family of $\kappa$-dense subsets of $X_{\bbb}$, hence we may
apply Theorem 4 to this $\bbb$ and $\dcal$ to get a family $\cbb$ of
$2^\kappa$ many $2$-partitions of ${\kappa}$ that satisfies
conditions (1) and (2) above.

The space that we need will be a further refinement of $X_{\cbb}$.
To obtain that, we next fix a 2-partition $\<I,J\>$ of the index set
$2^{\kappa}$ such that $|I|=|J|=2^{\kappa}$. For every unordered
pair $a\in \br I;2;$ we shall write $\ap=\max a$ and $\am=\min a$,
so that $a=\{\am,\ap\}$.

Let $\{j(a,m):a\in \br I;2;, m<{\omega}\}$ be pairwise distinct
elements of $J$. For any $a\in \br I;2; $ and $m<{\omega}$ we then
define the sets
\begin{displaymath}
E^0_{a,m}=D^0_{j(a,m)}\setm (D^0_{\am}\cap D^0_{\ap}) \text { and }
E^1_{a,m}={\kappa}\setm E^0_{a,m}.
\end{displaymath}
Clearly, then we have
\begin{displaymath}
E^1_{a,m}=D^1_{j(a,m)}\cup  (D^0_{\am}\cap D^0_{\ap}).
\end{displaymath}
In this way we obtained a new family
\begin{displaymath}
 \mathbb E=\left\{\<E^0_{a,m}, E^1_{a,m}\>:a\in \br
 I;2;,m<{\omega}\right\}
  \end{displaymath}
of 2-partitions of $\kappa$. We shall show that the space
$X_{\mathbb C\cup \mathbb E }$ satisfies all the requirements of
theorem \ref{tm:main2}.

\begin{claim}\label{cl:red}
For any finite function ${\eta}\in \fn(\br I;2;\times {\omega}, 2)$
and any ordinal ${\alpha}\in I$ there is a finite function
${\varphi}\in \fn(2^{\kappa},2)$ such that ${\alpha}\notin\dom
{\varphi}$ and $\eb[{\eta}]\supset \db[{\varphi}]$.
\end{claim}

\newcommand{\aaa}{a^*}

\begin{proof}[Proof of the Claim]
For each $a\in \br I;2;$ let us pick $a^*\in a$ with $a^* \ne
{\alpha}$. Then we have
\begin{multline}\notag
\eb[{\eta}]=
\bigcap_{{\eta}(a,m)=0}E^0_{a,m}\cap
\bigcap_{{\eta}(a,m)=1}E^1_{a,m}\sups
\\
\sups \bigcap_{{\eta}(a,m)=0} (D^0_{j(a,m)}\setm (D^0_{\am}\cap
D^0_{\ap})\cap \bigcap_{{\eta}(a,m)=1} D^1_{j(a,m)}\supset
\\
\sups\bigcap_{{\eta}(a,m)=0} (D^0_{j(a,m)}\cap D^1_{\aaa})\cap
\bigcap_{{\eta}(a,m)=1}
D^1_{j(a,m)}=\\
= \bigcap_{{\eta}(a,m)=0}D^1_{\aaa}\cap \bigcap_{\<a,m\>\in \dom
{\eta}}D^{{\eta}(a,m)}_{j(a,m)}.
\end{multline}
The expression in the last line above is, however, equal to
$\db[{\varphi}]$ for a suitable ${\varphi}\in \fn(2^{\kappa},2)$
because $j$ is an injective map of $[I] \times \omega$ into $J$ and
$\aaa\ne \alpha$ belongs to $I = \kappa \setminus J$ for all $a \in
[I]^2$.
\end{proof}

\begin{claim}
$\mathbb C\cup \mathbb E$ is ${\kappa}$-independent, hence
$\Delta(X_{\mathbb C\cup \mathbb E })={\kappa}$.
\end{claim}

\begin{proof}[Proof of the Claim]
Let ${\varepsilon}\in \fn (2^{\kappa},2)$ and ${\eta}\in \fn(\br
I;2;\times {\omega}, 2)$ be picked arbitrarily. By Claim
\ref{cl:red} there is ${\varphi}\in \fn(2^{\kappa},2)$ such that
$\eb[{\eta}]\supset \db[{\varphi}]$. Since $\db[{\varphi}]\in \dcal$
we have $|\cb[{\varepsilon}]\cap \db[{\varphi}]|={\kappa}$ because
$\cb$ satisfies condition (1). Consequently, we have
$|\cb[{\varepsilon}]\cap \eb[{\eta}]|={\kappa}$ as well.
\end{proof}

\begin{claim}
The family $\{D^0_{\alpha}:{\alpha}\in I\}$ witnesses that
$X_{\mathbb C\cup \mathbb E }$ is almost $2^{\kappa}$-resolvable.
\end{claim}

\begin{proof}[Proof of the Claim]
First we show that $D^0_{\alpha}$ is dense in $X_{\mathbb C\cup
\mathbb E }$ whenever ${\alpha}\in I$. So fix ${\alpha}\in I$,
moreover let ${\varepsilon}\in \fn (2^{\kappa},2)$  and ${\eta} \in
\fn(\br I;2;\times {\omega}, 2)$. By Claim \ref{cl:red} there is
${\varphi}\in \fn(2^{\kappa},2)$ such that ${\alpha}\notin\dom
{\varphi}$ and $\eb[{\eta}]\supset \db[{\varphi}]$. Since
${\alpha}\notin\dom {\varphi}$ we have $D^0_{\alpha}\cap
\db[{\varphi}] \in \dcal$. Hence, as $\cb$ has property (1),
$$\empt\ne (D^0_{\alpha}\cap \db[{\varphi}])\cap \cb[{\varepsilon}]\subs
D^0_{\alpha}\cap (\eb[{\eta}]\cap \cb[{\varepsilon}])$$ as well. So
$D^0_{\alpha}$ intersects every basic open subset of $X_{\cb\cup
\eb}$, i. e. $D^0_{\alpha}$ is dense in $X_{\cb\cup \eb}$.

Next we show that $D_{\alpha}\cap D_{\beta}$ is nowhere dense in the
space $X_{\cb\cup \eb}$ whenever $a = \{{\alpha},{\beta}\}\in \br
I;2;$. Indeed, let $\cb[{\varepsilon}]\cap \eb[{\eta}]$ be again a
basic open set with ${\varepsilon}\in \fn (2^{\kappa},2)$ and
${\eta}\in \fn (\br I;2;\times {\omega},2)$ and let us pick
$m<{\omega}$ such that $\<a,m\>\notin \dom {\eta}$. Then
$${\eta}'={\eta}\cup\{\<\<a,m\>,0\>\}
\in \fn(\br I;2;\times {\omega},2),$$ hence $\cb[{\varepsilon}]\cap
\eb[{\eta}']\subs \cb[{\varepsilon}]\cap \eb[{\eta}]$ is a
(non-empty) basic open set in the space $X_{\cb\cup \eb}$. Moreover,
$E^0_{a,m}=D^0_{j(a,m)}\setm (D^0_{\alpha}\cap D^0_{\beta})$ implies
\begin{multline}\notag
(D_{\alpha}\cap D_{\beta})\cap \mathbb C[{\varepsilon}] \cap \mathbb
E[{\eta}']\subs (D_{\alpha}\cap D_{\beta}) \cap (D^0_{j(a,m)}\setm
(D_{\alpha}\cap D_{\beta}))=\emptyset,
\end{multline}
consequently, $D_{\alpha}\cap D_{\beta}$ is not dense in
$\cb[{\varepsilon}]\cap \eb[{\eta}]$.
\end{proof}

Finally, the following simple claim will complete the proof of our
theorem.
\begin{claim}
The space $X_{\mathbb C }$ is $\omega_1$-irresolvable, that is, not
${\omega}_1$-resolvable.
\end{claim}
\begin{proof}[Proof of the Claim]
Assume that $\{F_{\zeta}:{\zeta}<{\omega}_1\}$ is a family of dense
subsets of $X_{\cb}$. By condition (2) the topology of $X_{\cb}$ is
$\dcal$-forced, so every $F_{\zeta}$ includes a $\dcal$-mosaic in
$X_{\cb}$, consequently for all $\zeta < \omega_1$ there are
${\varepsilon}_{\zeta}\in \fn(2^{\kappa},2)$ and $\phi_{\zeta}\in
\fn (2^{\kappa},2)$ such that $\db[\phi_{\zeta}] \cap
\cb[{\varepsilon}_{\zeta}]\subs F_{\zeta}$. By the well-known
$\Delta$-system lemma we may then find ${\zeta}<{\xi}<{\omega}_1$
such that ${\varepsilon}={\varepsilon}_{\zeta}\cup
{\varepsilon}_{\xi} \in \fn(2^{\kappa},2)$ and
$\phi=\phi_{\zeta}\cup \phi_{\xi}\in \fn(2^{\kappa},2)$. (Actually,
much more is true: there is an uncountable set $S \in
[\omega_1]^{\omega_1}$ such that the members of both
$\{{\varepsilon}_{\zeta} : \zeta \in S\}$ and $\{\phi_{\zeta} :
\zeta \in S\}$ are pairwise compatible.) But then we have
$$F_{\zeta}\cap F_{\xi}\supset \db[\phi_{\zeta}] \cap
\cb[{\varepsilon}_{\zeta}]\cap \db[\phi_{\xi}] \cap
\cb[{\varepsilon}_{\xi}]= \db[\phi]\cap \cb[\phi]\ne \empt.$$
\end{proof}
To conclude our proof, it suffices to recall the obvious fact that
if a topology on a set is $\lambda$-resolvable then so is any
coarser topology. Hence the $\omega_1$-irresolvability of
$X_{\mathbb C }$ implies that of $X_{\mathbb C\cup \mathbb E }$.
\end{proof}

Let us point out that as extraresolvability implies almost
$\omega$-resolvability that is equivalent to $\omega$-resolvability,
any counterexample to problem \ref{pr:ex} is automatically an
example of an $\omega$-resolvable but not maximally resolvable
space, hence it is a solution to the celebrated problem of Ceder and
Pearson from \cite{CP}. The first Tychonov ZFC examples of such
spaces were given in \cite{d-forced} and the spaces constructed in
theorem \ref{tm:main2} extend the supply of such examples.

\end{document}